\documentclass[11pt,reqno]{amsart}
\usepackage{}
\usepackage{mathrsfs}

\usepackage{amssymb,amsmath,epsfig}
\usepackage{amsfonts}
\usepackage{pst-plot}

\newcommand{\Z}{\mathbb{Z}}
\newcommand{\Q}{\mathbb{Q}}

\def\H{\mathbb{H}}

\theoremstyle{plain}
\newtheorem{thm}{Theorem}[section]

\newtheorem{prop}[thm]{Proposition}

\theoremstyle{definition}

\numberwithin{equation}{section}

\title{A short proof of an identity for cubic partitions }
\date{\today}
\author{Xinhua, Xiong}
\address{Department of Mathematics, China Three Gorges University, Yichang 443002,
P.R. China .}

\email{xinhuaxiong@ctgu.edu.cn}

\begin{document}

\begin{abstract}

In this note, we will give a short proof of an identity for  cubic partitions.

\end{abstract}
\maketitle

\section{Introduction}
\label{}
Let $p(n)$ denote the number of the unrestricted partitions of $n$, defined by $\sum_{n=0}^{\infty}q^n=
\prod_{n=0}^{\infty}\frac{1}{1-q^n}$. One of the celebrated results about $p(n)$ is the theorem which was proved by
Watson \cite{Watson38}: if $k\geq 1$, then for every nonnegative integer $n$
\begin{equation}\label{watson}
p(5^{k}n+r_k) \equiv 0 \pmod{5^k},
\end{equation}
where $r_k$ is the reciprocal modulo $5^k$ of $24$. Recently, the notion of cubic partitions of a natural number $n$ , named by
Kim \cite{Kim08}, was introduced by Chan \cite{Chan08a} in connection with Ramanujan's cubic continued fraction. By defintion, the generating function of
the number of cubic partitions of $n$ is
\begin{equation}
\sum_{n=0}^{\infty}a(n)q^n=
\prod_{n=0}^{\infty}\frac{1}{(1-q^n)(1-q^{2n})}.
\end{equation}
Chan \cite{Chan08a} from the Ramanujan's cubic continued fraction
$$
v(q):= \frac{q^{\frac{1}{3}}}{1}_+\frac{q+q^3}{1}_+\frac{q^2+q^4}{1}_+\dots \quad\quad |q| \le 1
$$
derived an elegant identity: let $x(q)= q^{-\frac{1}{3}}v(q)$, then
$$
\frac{1}{x(q)}-q^{\frac{1}{3}}-2q^{\frac{2}{3}}x(q)=\frac{(q^{\frac{1}{3}}; q^{\frac{1}{3}})_\infty(q^{\frac{2}{3}};q^{\frac{2}{3}})_\infty}{(q^3; q^3)_\infty(q^6; q^6)_\infty},
$$
where we set for $|q| \le1$, $(c;q)_\infty:=\prod_{n=0}^{\infty}(1-cq^k)$.
From this he obtained the generating function for $a(3n+2)$ \cite[Theorem 1]{Chan08a}:
\begin{thm}\label{1.1}
\begin{eqnarray}\label{1}
\sum_{n=0}^{\infty}a(3n+2)q^n=3\frac{(q^{3}; q^{3})_\infty^3(q^{6};q^{6})^3_\infty}{(q; q)^4_\infty(q^2; q^2)^4_\infty}.
\end{eqnarray}
\end{thm}
 Moreover he use this and the method of Hirschhorn and
Hunt derived the following congruences \cite[Theorem 1]{Chan08b} analogous to (\ref{watson}) :

\begin{thm}\label{1.2}
For every $k\ge 1$ and nonnegative integer $n$,
\begin{equation}\label{genconga}
a(3^{\alpha}n+c_{\alpha})\equiv  0\ ({\rm mod}\ 3^{\alpha+\delta(\alpha)}),
\end{equation}
where $c_{\alpha}$ is the reciprocal modulo $3^{\alpha}$ of $8$, and $\delta(\alpha)=1$
if $\alpha$ is even and $\delta(\alpha)=0$ otherwise.
\end{thm}
In this note, we will give a short proof of Theorem \ref{1.1} and outline the proof of Theorem \ref{1.2}.  The proof of Theorem \ref{1.1} by Chan used identities involved
Ramanujan's cubic continued fraction.  The proof of Theorem \ref{1.2} used the $H$ operators on formal power series  used by Hirschhorn and
Hunt \cite{{Hirschhorn}}. Our proofs will use  meromorphic  modular functions on $\Gamma_0(6)$ and $\Gamma_0(18)$.

\section{Preliminaries}
Let $\mathbb{H}:=\{z\in\mathbb{C}|{\rm Im}(z)>0\}$ denote the complex upper half plane, for a positive integer $N$, define the congruence subgroup $\Gamma_0(N)$ of $SL_2(\mathbb{Z})$ by
$
\Gamma_0(N):=\left\{\left(
                 \begin{array}{cc}
                   a & b \\
                   c & d \\
                 \end{array}
               \right)\Big|c\equiv 0\ ({\rm mod}\ N)
\right\}.
$
Let $\gamma=\left(
                 \begin{array}{cc}
                   a & b \\
                   c & d \\
                 \end{array}
               \right)\in SL_2(\mathbb{Z})$ act on the complex upper
half plane by the linear
fractional transformation
$
\gamma z:=\frac{az+b}{cz+d}.
$
Let $f(z)$ be a function on $\H$ which satisfies
$
f(\gamma z)=f(z),
$
if $f(z)$ is meromorphic on $\H$ and at  all the cusps of $\Gamma_0(N)$, then we call $f(z)$
a meromorphic modular function with respect to $\Gamma_0(N)$. The set of all such functions is denoted by
$\mathcal{M}_0(\Gamma_0(N))$.

Dedekind's eta function is defined by
$
\eta(z):=q^{\frac{1}{24}}\prod_{n=1}^\infty (1-q^n),
$
where $q=e^{2\pi iz}$ and ${\rm Im}(z)>0$.
A function $f(z)$ is called an eta-product if it can
be written in the form of
$
f(z)=\prod_{\delta |N}\eta^{r_{\delta}}(\delta z),
$
where $N$ is a natural number and  $r_{\delta}$ is an integer. The following
fact which is  due to Gordon, Hughes \cite{Gordon84} and Newman \cite{Newman} is useful  to verify  whether an eta-product is a
modular function.

\begin{prop}\label{prop2.1}
If $f(z)=\prod_{\delta|N}\eta^{r_{\delta}}(\delta
z)$ is an eta-product with
$
\frac{1}{2}\sum_{\delta|N}r_{\delta} =0
$
satisfies the  conditions:
$$
\sum_{\delta|N}\delta r_{\delta}\equiv 0 \ ({\rm mod}\ 24),\quad
\sum_{\delta|N}\frac{N}{\delta} r_{\delta}\equiv 0 \ ({\rm mod}\
24),\quad \prod_{\delta|N}\delta^{r_{\delta}}\in \Q^2,
$$
then $f(z)$ is in $\mathcal{M}_0(\Gamma_0(N))$.
\end{prop}

The following formula which is  due to Ligozat \cite{Ligozat75} gives the analytic orders of an eta-product at the
cusps of $\Gamma_0(N)$.

\begin{prop}\label{prop2.2}
Let $c,d$ and $N$ be positive integers with $d|N$ and $(c,d)=1$. If
$f(z)$ is an eta-product satisfying the conditions in Proposition
\ref{prop2.1} for $N$, then the order of vanishing of $f(z)$ at the
cusp $\frac{c}{d}$ is
$$\label{formula}
\frac{N}{24}\sum_{\delta
|N}\frac{(d,\delta)^2r_{\delta}}{(d,\frac{N}{d})d\delta}.
$$
\end{prop}

Let $p$ be a prime, and
$
f(q)= \sum_{n\ge n_0}^{\infty}a(n)q^n
$
be a formal power series, we define
$
U_pf(q)=\sum_{pn\ge n_0}a(pn)q^n.
$
\noindent  If $f(z)\in \mathcal{M}_0(\Gamma_0(N))$,
then $f(z)$ has an expansion at the point $i\infty$ of the form $f(z)=\sum_{n=n_0}^{\infty}a(n)q^n$
where $q=e^{2\pi iz}$ and ${\rm Im}(z)>0$. We call this expansion the Fourier series
of the $f(z)$. Moreover we define $U_pf(z)$ to be the result of applying $U_p$
to the Fourier series $f(z)$.

We use the results on the $U_3$-operator (we write $U$ for $U_3$ in the following) acting on modular functions on $\mathcal{M}_0(\Gamma_0(6))$ and $\mathcal{M}_0(\Gamma_0(18))$
stated by Gordon and Hughes \cite{Gordon84}. We know that $\Gamma_0(6)$ has $4$ cusps, represented by $0, \frac{1}{2}, \frac{1}{3}, \frac{1}{6}(=i\infty),$
$\Gamma_0(18)$ has $8$ cusps, represented by $0, \frac{1}{2}, \frac{1}{3}, \frac{2}{3}, \frac{1}{6},
\frac{5}{6}, \frac{1}{9}, \frac{1}{18} (= \infty).$ By Ligozat's
formula on the analytic orders of an eta-product, if $f(z)$ is in $\mathcal{M}_0(\Gamma_0(N))$, then $f(z)$ has the same order
at the cusps which have the same denominators. The order of $U_3(f(z))$ at a cusp $r$ of $\Gamma_0(6)$
is denoted by $ord_r U(f)$, and the order of $f(z)$ at a cusp of $s$ of $\Gamma_0(18)$ is denoted by $ord_s f$.

\begin{prop}\label{prop2.3}
Let $f(z)$ be an eta-product in $\mathcal{M}_0(\Gamma_0(18))$, then $U_3(f(z))$  is in $\mathcal{M}_0(\Gamma_0(6))$, and

\begin{eqnarray}
ord_0 U(f) &\geq& \mbox{min}\,(ord_0f, ord_{\frac{1}{3}}f), \quad
ord_{\frac{1}{2}} U(f) \geq\mbox{min}\,(ord_{\frac{1}{2}}f, ord_{\frac{1}{6}}f), \\
ord_{\frac{1}{3}} U(f) &\geq& \frac{1}{3}ord_{\frac{1}{9}}f, \quad\quad\quad\quad\quad\,\,\,
ord_{\frac{1}{6}} U(f) \geq \frac{1}{3}ord_{\frac{1}{18}}f.
\end{eqnarray}

Moreover, $U(f)$ has no poles on $\H$ except the cusps.
\end{prop}
\section{Proof of Theorem \ref{1.1} }
Let the eta-product
$$
F:=F(z)=\frac{\eta(9z)}{\eta(z)}\frac{\eta(18z)}{\eta(2z)},
$$
\noindent put $N=18$, we find $F(z)$ satisfies the conditions of Newman-Gordon-Hughes's theorem i.e. Proposition \ref{prop2.1},
so $F(z)$ is in $\mathcal{M}_0(\Gamma_0(18))$. We use Ligozat's formula to calculate the orders of $F(z)$
at the cusps $\frac{c}{d}$, for $ d=1,\,2,\,3,\,6,\,9,\,18$. We give the calculation of the case of  $d=1$ as an example :
\begin{eqnarray*}
ord_0F&=&\frac{18}{24\times (1, \frac{18}{1})} \sum_{\delta | 18}\frac{(1, \delta)^2}{\delta}r_{\delta}\\
&=&\frac{18}{24}\times \left(\frac{(1,9)^2}{9}\times 1 +\frac{(1, 18)^2}{18}\times 1 +\frac{(1, 1)^2}{1}\times (-1)+\frac{(1, 2)^2}{2}\times (-1) \right)\\
&=& -1.
\end{eqnarray*}
\noindent Similar calculations give
$$
ord_{\frac{1}{2}}F= -1,\quad ord_{\frac{1}{3}}F=0,\quad ord_{\frac{1}{6}}F=0,\quad ord_{\frac{1}{9}}F=1,\quad
ord_{\frac{1}{18}}F=1.
$$
\noindent By Proposition \ref{prop2.3}, the orders of $U(F)$ at the cusps of $\Gamma_0(6)$ satisfy
$$
ord_0U(F)\geq -1,\quad ord_{\frac{1}{2}}U(F) \geq -1,\quad
ord_{\frac{1}{3}}U(F)\geq 1,\quad ord_{\frac{1}{6}}U(F) \geq 1
$$
and $U(F)$ is holomorphic on $\H$. We note that the poles of $U(F)$ only appear at the cusps $0$ and $\frac{1}{2}$.
We define another eta-product
$$
A:=A(z)=\frac{\eta^4(3z)}{\eta^4(z)}\frac{\eta^4(6z)}{\eta^4(2z)}.
$$
By Proposition \ref{prop2.1}, we find that $A$ is in $\mathcal{M}_0(\Gamma_0(6))$. Ligozat's formula on the orders of an eta-product gives
$$
ord_0 A=-1,\quad ord_{\frac{1}{2}}A=-1,
ord_{\frac{1}{3}}A=1,\quad ord_{\frac{1}{6}}A=1
$$
and $A$ is holomorphic and non-zero elsewhere. Since the Riemann surface $(\H\cup \Q \cup{i\infty}) /\Gamma_0(6)$
has genus $0$, $\mathcal{M}_0(\Gamma_0(6))$ has one generator as a field. The orders of $A$ show that $U(F)=cA$ .
Since
\begin{eqnarray*}
F&=& q^3\prod_{n=1}^{\infty}\frac{(1-q^{9n})(1-q^{18n})}{(1-q^n)(1-q^{2n})}=q+q^2+3q^3 + 4q^4 +9q^5+12q^6+23q^7+31q^8+54q^9\dots,\\
A&= &q\prod_{n=1}^{\infty}\frac{(1-q^{3n})^4(1-q^{6n})^4}{(1-q^n)^4(1-q^{2n})^4}=q+4q^2+18q^3+52q^4+\dots.
\end{eqnarray*}
So $U(F)=3q+12q^2+54q^3+\dots $. The comparison of the coefficients of $U(f)$ and $A$ shows that $c=3$, so $U(F)=3A$.
On the other hand,
\begin{eqnarray*}
F&=& q^3\prod_{n=1}^{\infty}\frac{(1-q^{9n})(1-q^{18n})}{(1-q^n)(1-q^{2n})}
=\left(\sum_{n=1}^{\infty}a(n-1)q^n\right)\prod_{n=1}^{\infty}(1-q^{9n})(1-q^{18n})
\end{eqnarray*}
Apply $U$-operator again on both sides of above ,we have
\begin{equation}\label{U-1}
U(F)=3A =\left(\sum_{n=0}^{\infty}a(3n-1)q^n\right)\prod_{n=1}^{\infty}(1-q^{3n})(1-q^{6n}).
\end{equation}
 Put $$A=q\prod_{n=1}^{\infty}\frac{(1-q^{3n})^4(1-q^{6n})^4}{(1-q^n)^4(1-q^{2n})^4}$$
into above, we obtain the identity :
\begin{eqnarray}
\sum_{n=0}^{\infty}a(3n+2)q^n=3\frac{(q^{3}; q^{3})_\infty^3(q^{6};q^{6})^3_\infty}{(q; q)^4_\infty(q^2; q^2)^4_\infty}.
\end{eqnarray}
Which is Theorem \ref{1.1}.

We outline the proof of Theorem \ref{1.2}.  The ideal is similar to the paper \cite{Xiong}. Firstly we apply the proposition \ref{prop2.3} to $A^i $ and $FA^i$ for $i\geq 1$, we can express $A^i $ (resp. $FA^i$) as a polynomial in $A$ of degree at most $3i$ (resp. $3i+1$ ). This is the part corresponding to the Proposition 1 and Proposition 2 in \cite{Chan08b}. Next we use the initial values of $A^i $ to calculate the elementary symmetric functions $\sigma_i$ ($i=1,2,3$) of $U(A(\frac{z+t}{3}))\, (t=0,\,1,\,2)$ which are polynomials in $A$ with integers as coefficients. Then by the Newton recurrence for
power sums, we get for all $i\geq 3$
$$
U(A^i)= \sigma_1 U(A^{i-1})-\sigma_2 U(A^{i-2})+\sigma_3 U(A^{i-3})
$$
Hence for $i\geq 1$, $U(A^i) \in \Z[A]$. Moreover $U(FA^i)$ satisfies the same recurrence as $U(A^i)$ also $U(FA^i)$ is in $\Z[A]$ and By induction we obtain the lower bounds of  $3$-adic orders of these coefficients.  The last step  is almost the same as the Proposition $3$ and the Theorem $4$.


\end{document}